\documentclass[12pt,a4paper,reqno]{amsart} %Fr o m detta
\usepackage{amssymb}
\usepackage{latexsym}
\usepackage{amsmath}
\usepackage{mathrsfs}
\usepackage{amsthm}
\usepackage{calc}                         %Detta behövs för att fixa
\usepackage[english]{babel}
\usepackage{cite}
\usepackage{xcolor,soul}
\newcommand{\abs}[1]{\left|#1\right|}

\newcommand{\scal}[2]{\langle #1,#2\rangle}

\newcommand{\prn}[1]{\left(#1\right)}

\newcommand{\rr}[1]{\mathbf R^{#1}}

\newcommand{\nn}[1]{\mathbf N^{#1}}

\newcommand{\nm}[2]{\Vert #1\Vert _{#2}}

\newcommand{\op}{\operatorname{Op}}

\newcommand{\set}[1]{\left\{ \, #1\,\right \} }
\newcommand{\sets}[2]{\{ \, #1\, ;\, #2\, \} }

\newcommand{\fy}{\varphi}

\newcommand{\cdo}{\, \cdot \, }
\newcommand{\loc}{\operatorname{loc}}

\newcommand{\vrum}{\vspace{0.1cm}}
\newcommand{\GL}{\mathbf{M}}

\newcommand{\maclS}{\mathcal S}

\newcommand{\mascB}{\mathscr B}

\newcommand{\mascF}{\mathscr F}
\newcommand{\mascP}{\mathscr P}
\newcommand{\mascS}{\mathscr S}

\numberwithin{equation}{section}          %Detta gör att man får
\newtheorem{thm}{Theorem}
\numberwithin{thm}{section}

\newcommand{\rubrik}{}
\newtheorem{prop}[thm]{Proposition}

\newtheorem{lemma}[thm]{Lemma}

\theoremstyle{definition}

\newtheorem{defn}[thm]{Definition}
\newtheorem{example}[thm]{Example}

\theoremstyle{remark}

\newtheorem{rem}[thm]{Remark}              %T o m hit är bara allmän
\author{Ahmed Abdeljawad}

\address{Department of Mathematics,
University of Turin, Italy}

\email{ahmed.abdeljawad@unito.it}

\author{Joachim Toft}

\address{Department of Mathematics,
Linn{\ae}us University, Sweden}

\email{joachim.toft@lnu.se}

\title{\textbf{Anisotropic Gevrey-H\"ormander pseudo-differential operators
on modulation spaces}}

\begin{document}



\begin{abstract}
We show continuity properties for the pseudo-differential operator
$\op (a)$ from $M(\omega _0\omega ,\mascB )$
to $M(\omega ,\mascB )$,
for fixed $s,\sigma\ge 1$,
$\omega ,\omega _0\in \mascP _{s,\sigma}^0$
($\omega ,\omega _0\in \mascP _{s,\sigma}$), $a\in
\Gamma ^{\sigma,s}_{(\omega _0)}$
($a\in \Gamma ^{\sigma,s;0}_{(\omega _0)}$)
, and $\mascB$ is an invariant
Banach function space.
\end{abstract}

\maketitle

%%%%%%%%%%%%%%%%%%%%%%%%%%%%%%%%%
\section{Introduction}\label{sec0}
%%%%%%%%%%%%%%%%%%%%%%%%%%%%%%%%%

\par

In the paper, we consider pseudo-differential operators, where
the symbols are of infinite orders and possess suitable Gevrey regularities
and which are allowed to grow sub-exponentially together with all
their derivatives. Our main purpose is to extend boundedness results,
in \cite{To14}, of the
pseudo-differential operators when acting on modulation spaces.

\par

More specific, the symbols should satisfy conditions of the form
\begin{equation}\label{Eq:SymbCondIntr}
|\partial _x^\alpha \partial _\xi ^\beta a(x,\xi )|
\lesssim h^{|\alpha +\beta |}\alpha !^\sigma \beta !^s\omega _0(x,\xi ),
\end{equation}
where $\omega _0$ should be a moderate weight on $\rr {2d}$
and satisfy boundedness conditions like
\begin{equation}\label{Eq:WeightCondIntr}
\omega _0(x,\xi )\lesssim e^{r(|x|^{\frac 1s}+|\xi |^{\frac 1\sigma})}.
\end{equation}
For such symbols $a$ we prove that corresponding pseudo-differential
operators $\op (a)$ is continuous from the modulation space
$M(\omega _0\omega ,\mascB )$ to $M(\omega ,\mascB )$.
(See Section \ref{sec1} for notations.)

\par

Similar investigations were performed in \cite{To25} in the case
$s=\sigma$ (i.{\,}e. the isotropic case).
Therefore, the results in the current paper are more general
in the sense of the anisotropicity of the considered symbol classes.
Moreover, we use different techniques compared to \cite{To25}.

\par

We also remark that several ideas arise in \cite{To14}, where similar
investigations were performed after the conditions \eqref{Eq:SymbCondIntr}
and \eqref{Eq:WeightCondIntr} are replaced by
$$
|\partial _x^\alpha \partial _\xi ^\beta a(x,\xi )| \lesssim \omega _0(x,\xi )
$$
and
$$
\omega _0(x,\xi )\lesssim (1+|x|+|\xi |)^N,
$$
respectively, for some $N\ge 0$.

\par

In \cite{Fe4}, H. Feichtinger introduced the modulation spaces
to measure the
time-frequency concentration of a function or distribution
on the time-frequency
space or the phase space $\rr {2d}$.
Nowadays they become popular among mathematicians
and engineers since their numerous applications in
signal processing \cite{Fe8,Fe9}, pseudo-differential and Fourier integral operators
\cite{CoGaTo,CoTo,CoGrNiRo1,PT2,PT3,Ta,Te1,Te2,To14,To18,To20,To22,To24,To25}
and quantum mechanics \cite{CoGrNiRo,dGo}.

\par

%A pseudo-differential operator in $\rr d$ with symbol
%$a\in L^1(\rr d)$ is defined by the formula
%%%
%\begin{equation}\label{KN}
%\op (a)f(x)=(2\pi)^{-\frac d2}\int_{\rr d}
%a(x,\xi) \widehat f(\xi)e^{i\scal  x\xi}
%\,d\xi,\qquad
%f\in\mathscr S(\rr d).
%\end{equation}
%%%
%%where
%%$\widehat{f}(x)=\mathscr F{f}(x)
%%=
%%(2\pi)^{\frac d2}\int_{\rr d} e^{- i \scal x\xi} f(x)\,dx$
%%is the Fourier transform of $f$.
%The definition extends to any $a\in \mathscr S '(\rr d)$, and
%then $\op (a)f$ is well-defined as a temperate distribution.
%
%
%\par
%
%Pseudo-differential operators arise in different situations, e.{\,}g.
%partial differential equations, quantum mechanics and engineering.
%In the theory of partial differential equations, they were introduced
%in \cite{hormander65,kohn-nirenberg65}.
%Thereafter several symbol classes have been considered,
%depending of the applications of the  partial differential equations.
%In particular, a well deep analysis of such operators have
%been carried on for $\SG$ classes
%$S^{m,\mu}_{1,1}$, $m,\mu\in\re$, of smooth functions $a(x,\xi)$
%satisfying the estimates
%$|\partial^\alpha_x\partial^\beta_\xi
%a(x,\xi)|\leq
%C_{\alpha,\beta}\langle x \rangle^{m-|\alpha|}
%\langle \xi \rangle ^{\mu -|\beta|} $.

\par

The paper is organized as follows. In Section \ref{sec1}
we give the main definition and properties of Gelfand-Shilov and
modulation spaces and we recall some essential results.
In Section \ref{sec2} we state our main results on
the continuity with anisotropic settings.

\par

%%%%%%%%%%%%%%%%%%%%%%%%%%%%%%%%%
\section{Preliminaries}\label{sec1}
%%%%%%%%%%%%%%%%%%%%%%%%%%%%%%%%%

\par

In the current section we review basic properties for modulation
spaces and other related spaces. More details and proofs can be found in
\cite{Fe2,Fe3,Fe4,FG1,FG2,FG4,GaSa,Gc2,To20} .

\par

\subsection{Weight functions}\label{subsec1.1}

\par

A function $\omega$ on $\rr d$ is called a \emph{weight} or
\emph{weight function},
if $\omega ,1/\omega \in L^\infty _{\loc} (\rr d)$
are positive everywhere.
The weight $\omega$ on $\rr d$ is called $v$-moderate
for some weight $v$ on $\rr d$, if
\begin{equation}\label{vModerate}
\omega (x+y)\lesssim \omega (x)v(y),\quad x,y\in \rr d.
\end{equation}
If $v$ is even and satisfies \eqref{vModerate} with $\omega =v$,
then $v$ is called submultiplicative. 

\par

Let $s,\sigma>0$. Then we let $\mathscr P_E(\rr d)$ be the
set of all moderate weights on $\rr d$, $\mathscr P_s(\rr  d)$
($\mathscr P_s^0(\rr  d)$) be the set of all $\omega\in
\mathscr P_E(\rr d)$ such that 
$$
\omega (x+y)\lesssim \omega (x)e^{r|y|^{\frac 1s}}, \quad x,y\in \rr d,
$$
for some $r>0$ (for every $r>0$), and  $\mathscr P_{s,\sigma}(\rr {2d})$
($\mathscr P_{s,\sigma}^0(\rr {2d})$) be the set of all
$\omega\in \mathscr P_E(\rr {2d})$ such that
\begin{equation} \label{estomega}
\omega (x+y,\xi+\eta)\lesssim \omega (x,\xi)
e^{r(|y|^{\frac 1s}+|\eta|^{\frac 1\sigma})}, \quad x,y,\xi,\eta\in \rr d,
\end{equation}
for some $r>0$ (for every $r>0$).

\par

%
%For any $s>0$, let
%$\mascP _s(\rr d)$ ($\mascP _s^0(\rr d)$) be the set of all
%weights $\omega$ on $\rr d$ such that
%$$
%e^{-r|x|^{\frac 1s}}\lesssim \omega (x)\lesssim e^{r|x|^{\frac 1s}}
%$$
%for some $r>0$ (for every $r>0$). In similar ways, if $s,\sigma >0$, then
%$\mascP _{s,\sigma}(\rr {2d})$ ($\mascP _{s,\sigma}^0(\rr {2d})$) consists of all
%submultiplicative weight functions $\omega$ on $\rr {2d}$ such that
%%$$
%%e^{-r(|x|^{\frac 1s}+|\xi |^{\frac 1\sigma})}\lesssim \omega (x,\xi )
%%\lesssim e^{r(|x|^{\frac 1s}+|\xi |^{\frac 1\sigma})}
%%$$
%%for some $r>0$ (for every $r>0$).
%%In particular, if $\omega \in \mascP _{s,\sigma}(\rr {2d})$
%($\mascP _{s,\sigma}^0(\rr {2d})$), then 
%\begin{equation} \label{estomega}
%\omega (x+y,\xi +\eta ) \lesssim \omega (x,\xi )e^{r(|y|^{\frac 1s}
%	+|\eta |^{\frac 1\sigma})}, \quad  x,y,\xi ,\eta  \in \rr d,
%\end{equation}
%%%
%for some $r>0$ (for every $r>0$).
%
%\par

The following result shows that for any weight in $\mascP _E$,
there are equivalent weights that 
satisfy strong Gevrey regularity.

\par

\begin{prop}\label{Prop:EquivWeights}
Let $\omega \in \mascP _E(\rr {2d})$ and $s,\sigma> 0$.
Then there exists a weight
$\omega _0\in \mascP _E(\rr {2d})\cap C^\infty (\rr {2d})$
such that the following is true:	
\begin{enumerate}	
\item $\omega _0\asymp \omega $;

\vrum

\item for every $h>0$,
\begin{equation*}
|\partial _x ^{\alpha}\partial _\xi ^\beta
\omega _0(x, \xi)|
\lesssim 
h^{|\alpha +\beta|}\alpha !^\sigma
\beta!^s \omega _0(x,\xi) 
\asymp 
h^{|\alpha +\beta|}\alpha !^\sigma
\beta!^s\omega (x,\xi).
\end{equation*} 
\end{enumerate}
\end{prop}

\par

Proposition \ref{Prop:EquivWeights} is equivalent to \cite[Proposition 1.6]{AbCoTo}.
In fact, by Proposition \cite[Proposition 1.6]{AbCoTo} we have that
Proposition \ref{Prop:EquivWeights} holds with $s=\sigma$.
Hence, Proposition \ref{Prop:EquivWeights} implies \cite[Proposition 1.6]{AbCoTo}.
On the other hand, let $s_0=\min (s,\sigma)$. Then
\cite[Proposition 1.6]{AbCoTo} implies that there is a weight function
$\omega_0\asymp \omega$ satisfying
\begin{align*}
|\partial _x ^{\alpha}\partial _\xi ^\beta
\omega _0(x, \xi)|
&\lesssim 
h^{|\alpha +\beta|}(\alpha !
\beta!)^{s_0} \omega _0(x,\xi) 
\\[1ex]
&\lesssim
h^{|\alpha +\beta|}\alpha !^\sigma
\beta!^s \omega _0(x,\xi),
\end{align*}
giving Proposition \ref{Prop:EquivWeights}.

\par

\subsection{Gelfand-Shilov spaces}\label{subsec1.2}

\par

%We start by recalling some facts about Gelfand-Shilov spaces.
Let $0<h,s,\sigma \in \mathbf R$ be fixed. Then
$\maclS _{s;h}^{\sigma}(\rr d)$ is
the Banach space of all $f\in C^\infty (\rr d)$ such that
\begin{equation}\label{gfseminorm}
\nm f{\maclS _{s;h}^{\sigma}}\equiv \sup_{\alpha ,\beta \in
\mathbf N^d} \sup_{x \in \rr d} \frac {|x^\alpha \partial ^\beta
f(x)|}{h^{|\alpha | + |\beta |}\alpha !^s\, \beta !^\sigma}<\infty,
\end{equation}
endowed with the norm \eqref{gfseminorm}.

\par

The \emph{Gelfand-Shilov spaces} $\maclS _{s}^{\sigma}(\rr d)$ and
$\Sigma _{s}^{\sigma}(\rr d)$ are defined as the inductive and projective 
limits respectively of $\maclS _{s;h}^{\sigma}(\rr d)$. This implies that
\begin{equation}\label{GSspacecond1}
\maclS _{s}^{\sigma}(\rr d) = \bigcup _{h>0}\maclS _{s;h}^{\sigma}(\rr d)
\quad \text{and}\quad \Sigma _{s}^{\sigma}(\rr d) =\bigcap _{h>0}
\maclS _{s;h}^{\sigma}(\rr d),
\end{equation}
and that the topology for $\maclS _{s}^{\sigma}(\rr d)$ is the strongest
possible one such that the inclusion map from $\maclS _{s;h}^{\sigma}
(\rr d)$ to $\maclS _{s}^{\sigma}(\rr d)$ is continuous, for every choice 
of $h>0$. The space $\Sigma _{s}^{\sigma}(\rr d)$ is a Fr{\'e}chet space
with seminorms $\nm \cdo {\maclS _{s;h}^{\sigma}}$, $h>0$. Moreover,
$\Sigma _{s}^{\sigma}(\rr d)\neq \{ 0\}$, if and only
if $s+\sigma \ge 1$ and $(s,\sigma )\neq
(\frac 12,\frac 12)$, and
$\maclS _{s}^{\sigma}(\rr d)\neq \{ 0\}$, if and only
if $s+\sigma \ge 1$.

%\par
%
%The spaces $\maclS _{s}^{\sigma}(\rr d)$ and $\Sigma _{s}^{\sigma}(\rr d)$
%can be characterized 
%also in terms of the exponential decay of their elements,
%namely $f \in \maclS _{s}^{\sigma}(\rr d)$ 
%(respectively $f \in \Sigma _{s}^{\sigma}(\rr d)$) if and only if 
%$$
%|\partial^\alpha f(x)| \lesssim \ep^{|\alpha|} (\alpha!)^\sigma
%e^{-h|x|^{\frac 1s}}
%$$
%for some $h>0, \ep>0$ (respectively for every $h>0, \ep>0$). 
%Moreover we recall that for $s <1$ the elements of
%$\maclS _{s}^{\sigma}(\rr d)$ admit entire 
%extensions to $\mathbf{C}^d$ satisfying suitable exponential
%bounds, cf. \cite{GS} for details.

\medspace

The \emph{Gelfand-Shilov distribution spaces} $(\maclS _{s}^{\sigma})'(\rr d)$
and $(\Sigma _{s}^{\sigma})'(\rr d)$
are the projective and inductive limit
respectively of $(\maclS _{s;h}^{\sigma})'(\rr d)$.
In \cite{Pi2} it is proved that $(\maclS _{s}^{\sigma})'(\rr d)$
is the dual of $\maclS _s^\sigma(\rr d)$, and $(\Sigma _{s}^{\sigma})'(\rr d)$
is the dual of $\Sigma _{s}^{\sigma}(\rr d)$ (also in topological sense).
%
% This means that
%%%
%\begin{equation}\tag*{(\ref{GSspacecond1})$'$}
%(\maclS _{s}^{\sigma})'(\rr d) = \bigcap _{h>0}(\maclS _{s;h}^{\sigma})'(\rr d)\quad
%\text{and}\quad (\Sigma _{s}^{\sigma})'(\rr d) =\bigcup _{h>0}(\maclS _{s;h}^{\sigma})'(\rr d).
%\end{equation}
%%
%We remark that in \cite{Pi2} it is proved that $(\maclS _{s}^{\sigma})'(\rr d)$
%is the dual of $\maclS _s^\sigma(\rr d)$, and $(\Sigma _{s}^{\sigma})'(\rr d)$
%is the dual of $\Sigma _{s}^{\sigma}(\rr d)$ (also in topological sense).
%
%\par
%
%For every $s,\sigma >0$ we have
%%%
%\begin{equation}\label{Eq:GSEmbeddings}
%\Sigma _s^\sigma (\rr d)
%\hookrightarrow
%\maclS _s^\sigma (\rr d)
%\hookrightarrow
%\Sigma _{s+\ep}^{\sigma +\ep}(\rr d)
%\hookrightarrow
%\mascS (\rr d)
%\end{equation}
%%%
%for every $\ep >0$. If $s+\sigma \ge 1$, then
%the last two inclusions in \eqref{Eq:GSEmbeddings} are dense,
%and if in addition $(s,\sigma )\neq (\frac 12,\frac 12)$, then the
%first inclusion in \eqref{Eq:GSEmbeddings} is dense.
%
%\par
%
%From these properties it follows that $\mascS '(\rr d)\hookrightarrow
%(\maclS _s^\sigma)'(\rr d)$ when $s+\sigma \ge 1$, and if in addition
%$(s,\sigma )\neq (\frac 12,\frac 12)$, then $(\maclS _s^\sigma)'(\rr d)
%\hookrightarrow (\Sigma _s^\sigma)'(\rr d)$.

%\par
%
%The Gelfand-Shilov spaces possess several convenient mapping
%properties. For example they are invariant under
%translations, dilations, and to some extent tensor products
%and (partial) Fourier transformations.

\par

The Fourier transform $\mathscr F$ is the linear and continuous
map on $\mascS (\rr d)$,
given by the formula
$$
(\mathscr Ff)(\xi )= \widehat f(\xi ) \equiv (2\pi )^{-\frac d2}\int _{\rr
	{d}} f(x)e^{-i\scal  x\xi }\, dx
$$
when $f\in \mascS (\rr d)$. Here $\scal \cdo \cdo$ denotes the usual
scalar product on $\rr d$. 
The Fourier transform extends  uniquely to homeomorphisms
from $(\maclS _{s}^{\sigma})'(\rr d)$ to $(\maclS _{\sigma}^{s})'(\rr d)$,
and from  $(\Sigma _{s}^{\sigma})'(\rr d)$ to $(\Sigma _{\sigma}^{s})'(\rr d)$.
Furthermore, it restricts to homeomorphisms from
$\maclS _{s}^{\sigma}(\rr d)$ to $\maclS _{\sigma}^{s}(\rr d)$,
and from  $\Sigma _{s}^{\sigma}(\rr d)$ to $\Sigma _{\sigma}^{s}(\rr d)$.

\par

\medspace

Some considerations later on involve a broader family of
Gelfand-Shilov spaces. More precisely, for $s_j,\sigma _j\in \mathbf R_+$,
$j=1,2$, the Gelfand-Shilov spaces $\maclS _{s _1,s_2}^{\sigma _1,\sigma _2}(\rr {d_1+d_2})$ and
$\Sigma _{s _1,s_2}^{\sigma _1,\sigma _2}(\rr {d_1+d_2})$ consist of all functions
$F\in C^\infty (\rr {d_1+d_2})$ such that
\begin{equation}\label{GSExtCond}
|x_1^{\alpha _1}x_2^{\alpha _2}\partial _{x_1}^{\beta _1}
\partial _{x_2}^{\beta _2}F(x_1,x_2)| \lesssim
h^{|\alpha _1+\alpha _2+\beta _1+\beta _2|}
\alpha _1!^{s_1}\alpha _2!^{s_2}\beta _1!^{\sigma _1}\beta _2!^{\sigma _2}
\end{equation}
for some $h>0$ respective for every $h>0$. The  topologies, and the duals
\begin{alignat*}{3}
&(\maclS _{s _1,s_2}^{\sigma _1,\sigma _2})'(\rr {d_1+d_2}) &
&\quad \text{and} \quad &
&(\Sigma _{s _1,s_2}^{\sigma _1,\sigma _2})'(\rr {d_1+d_2})
\intertext{of}
&\maclS _{s _1,s_2}^{\sigma _1,\sigma _2}(\rr {d_1+d_2}) &
&\quad \text{and} \quad &
&\Sigma _{s _1,s_2}^{\sigma _1,\sigma _2}(\rr {d_1+d_2}),
\end{alignat*}
respectively, and their topologies
are defined in analogous ways as for the spaces $\maclS _s^\sigma (\rr d)$
and $\Sigma _s^\sigma (\rr d)$ above.

\par

The following proposition explains mapping properties of partial
Fourier transforms on Gelfand-Shilov spaces, and follows by similar
arguments as in analogous situations in
\cite{GS}. The proof is therefore omitted. Here, $\mascF _1F$
and $\mascF _2F$ are the partial
Fourier transforms of $F(x_1,x_2)$ with respect to
$x_1\in \rr {d_1}$ and $x_2\in \rr {d_2}$,
respectively.

\par

\begin{prop}\label{propBroadGSSpaceChar}
Let $s_j,\sigma _j >0$, $j=1,2$.
Then the following is true:
\begin{enumerate}	
\item the mappings $\mascF _1$ and $\mascF _2$ on $\mascS (\rr {d_1+d_2})$
restrict to homeomorphisms
\begin{align*}
\mascF _1 \, &: \, \maclS _{s _1,s_2}^{\sigma _1,\sigma _2}(\rr {d_1+d_2}) \to
\maclS _{\sigma _1,s_2}^{s_1,\sigma _2}(\rr {d_1+d_2})
\intertext{and}
\mascF _2 \, &: \, \maclS _{s _1,s_2}^{\sigma _1,\sigma _2}(\rr {d_1+d_2}) \to
\maclS _{s _1,\sigma _2}^{\sigma _1,s_2}(\rr {d_1+d_2})
\text ;
\end{align*}

\vrum

\item the mappings $\mascF _1$ and $\mascF _2$ on
$\mascS (\rr {d_1+d_2})$ are uniquely extendable to
homeomorphisms
\begin{align*}
\mascF _1 \, &: \, (\maclS _{s _1,s_2}^{\sigma _1,\sigma _2})'(\rr {d_1+d_2}) \to
(\maclS _{\sigma _1,s_2}^{s_1,\sigma _2})'(\rr {d_1+d_2})
\intertext{and}
\mascF _2 \, &: \, (\maclS _{s _1,s_2}^{\sigma _1,\sigma _2})'(\rr {d_1+d_2}) \to
(\maclS _{s _1,\sigma _2}^{\sigma _1,s_2})'(\rr {d_1+d_2}).
\end{align*}
\end{enumerate}
	
\par
	
The same holds true if the $\maclS  _{s _1,s_2}^{\sigma _1,\sigma _2}$-spaces and
their duals are replaced by
corresponding $\Sigma  _{s _1,s_2}^{\sigma _1,\sigma _2}$-spaces and their duals.
\end{prop}

\par

The next two results follow from \cite{ChuChuKim}. The proofs are therefore omitted.

\begin{prop}
Let $s_j,\sigma _j> 0$, $j=1,2$. Then the following
conditions are equivalent.
\begin{enumerate}
\item $F\in \maclS _{s _1,s_2}^{\sigma _1,\sigma _2}(\rr {d_1+d_2})$\quad
($F\in \Sigma _{s _1,s_2}^{\sigma _1,\sigma _2}(\rr {d_1+d_2})$);

\vrum

\item for some $h>0$ (for every $h>0$) it holds
\begin{equation*}
\displaystyle{|F(x_1,x_2)|\lesssim e^{-h(|x_1|^{\frac 1{s_1}} + |x_2|^{\frac 1{s_2}} )}}
\quad \text{and}\quad 
\displaystyle{|\widehat F(\xi _1,\xi _2)|\lesssim
	e^{-h(|\xi _1|^{\frac 1{\sigma _1}} + |\xi _2|^{\frac 1{\sigma _2}} )}}.
\end{equation*}
\end{enumerate}
\end{prop}

\par

We notice that if
$s_j+\sigma _j<1$ for some $j=1,2$, then
$\maclS _{s _1,s_2}^{\sigma _1,\sigma _2}(\rr {d_1+d_2})$
and $\Sigma _{s _1,s_2}^{\sigma _1,\sigma _2}(\rr {d_1+d_2})$
are equal to the trivial space $\{ 0\}$.
Likewise, if $s_j=\sigma _j=\frac 12$ for some $j=1,2$, then
$\Sigma _{s _1,s_2}^{\sigma _1,\sigma _2}(\rr {d_1+d_2}) = \{ 0\}$.

\par

\subsection{Short time Fourier transform and Gelfand-Shilov spaces}

\par

We recall here some basic facts about
the short-time Fourier transform and weights.

\par

Let $\phi \in \maclS _s^\sigma (\rr d)\setminus 0$ be fixed. Then the short-time
Fourier transform of $f\in (\maclS _s^\sigma )'(\rr d)$ is given by
$$
(V_\phi f)(x,\xi ) = (2\pi )^{-\frac d2}(f,\phi (\cdo -x)
e^{i\scal \cdo \xi})_{L^2}.
$$
Here $(\cdo ,\cdo )_{L^2}$ is the unique extension of the $L^2$-form on
$\maclS _s^\sigma (\rr d)$ to a continuous sesqui-linear form on $(\maclS
_s^\sigma )'(\rr d)\times \maclS _s^\sigma (\rr d)$. In the case
$f\in L^p(\rr d)$, for some $p\in [1,\infty]$, then $V_\phi f$ is given by
$$
V_\phi f(x,\xi ) \equiv (2\pi )^{-\frac d2}\int _{\rr d}f(y)\overline{\phi (y-x)}
e^{-i\scal y\xi}\, dy .
$$

\par

The following characterizations of the
$\maclS _{s_1,s_2}^{\sigma _1,\sigma _2}(\rr {d_1+d_2})$,
$\Sigma _{s_1,s_2}^{\sigma _1,\sigma _2}(\rr {d_1+d_2})$
and their duals
follow by similar arguments as in the proofs of
Propositions 2.1 and 2.2 in \cite{To22}. The details are left
for the reader.

\par

\begin{prop}\label{Prop:STFTGelfand2}
Let $s_j,\sigma _j>0$ be such that $s_j+\sigma _j\ge 1$, $j=1,2$, 
$s_0\le s$ and $\sigma_0\le \sigma$. Also let
$\phi \in \mathcal S_{s_1,s_2}^{\sigma _1,\sigma _2}(\rr {d_1+d_2}
\setminus 0$)
($\phi \in \Sigma _{s_1,s_2}^{\sigma _1,\sigma _2}(\rr {d_1+d_2}
\setminus 0$) and let $f$ be a Gelfand-Shilov distribution on
$\rr {d_1+d_2}$.
Then $f\in  \maclS _{s_1,s_2}^{\sigma _1,\sigma _2}(\rr {d_1+d_2})$
($f\in  \Sigma _{s_1,s_2}^{\sigma _1,\sigma _2}(\rr {d_1+d_2})$),
if and only if
\begin{equation}\label{stftexpest2}
|V_\phi f(x_1,x_2,\xi _1,\xi _2)|
\lesssim
e^{-r (|x_1|^{\frac 1{s_1}} + |x_2|^{\frac 1{s_2}}
	+|\xi _1|^{\frac 1{\sigma _1}} +|\xi _2|^{\frac 1{\sigma _2}} )},
\end{equation}
holds for some $r > 0$ (holds for every $r > 0$).
%
%\vrum
%
% if in addition
%$\phi \in \Sigma _{s_1,s_2}^{\sigma _1,\sigma _2}(\rr {d_1+d_2})
%\setminus 0$, then 
%$f\in  \Sigma _{s_1,s_2}^{\sigma _1,\sigma _2}(\rr {d_1+d_2})$
%if and only if
%%%
%\begin{equation}\label{stftexpest2A}
%|V_\phi f(x_1,x_2,\xi _1,\xi _2)|
%\lesssim
%e^{-r (|x_1|^{\frac 1{s_1}} + |x_2|^{\frac 1{s_2}}
%	+|\xi _1|^{\frac 1{\sigma _1}} +|\xi _2|^{\frac 1{\sigma _2}} )}
%\end{equation}
%%%
%holds for every $r > 0$.
%
\end{prop}

\par

A proof of Proposition  \ref{Prop:STFTGelfand2} can be found in
e.{\,}g. \cite{GZ} (cf. \cite[Theorem 2.7]{GZ}). The
corresponding result for Gelfand-Shilov distributions
is the following improvement of \cite[Theorem 2.5]{To18}.

\par

\begin{prop}\label{Prop:STFTGelfand2Dist}
Let $s_j,\sigma _j>0$ be such that $s_j+\sigma _j\ge 1$, $j=1,2$, 
$s_0\le s$ and $t_0\le t$. Also let
$\phi \in \mathcal S_{s_1,s_2}^{\sigma _1,\sigma _2}(\rr {d_1+d_2})
\setminus 0$ and let $f$ be a Gelfand-Shilov distribution on
$\rr {d_1+d_2}$.
Then the following is true:
	\begin{enumerate}
\item $f\in  (\maclS _{s_1,s_2}^{\sigma _1,\sigma _2})'(\rr {d_1+d_2})$,
if and only if
\begin{equation}\label{stftexpest2Dist}
|V_\phi f(x_1,x_2,\xi _1,\xi _2)|
\lesssim
e^{r (|x_1|^{\frac 1{s_1}} + |x_2|^{\frac 1{s_2}}
	+|\xi _1|^{\frac 1{\sigma _1}} +|\xi _2|^{\frac 1{\sigma _2}} )}
\end{equation}
holds for every $r > 0$;

\vrum

\item if in addition
$\phi \in \Sigma _{s_1,s_2}^{\sigma _1,\sigma _2}(\rr {d_1+d_2})
\setminus 0$, then 
$f\in  (\Sigma _{s_1,s_2}^{\sigma _1,\sigma _2})'(\rr {d_1+d_2})$,
if and only if
\begin{equation}\label{stftexpest2DistA}
|V_\phi f(x_1,x_2,\xi _1,\xi _2)|
\lesssim
e^{r (|x_1|^{\frac 1{s_1}} + |x_2|^{\frac 1{s_2}}
	+|\xi _1|^{\frac 1{\sigma _1}} +|\xi _2|^{\frac 1{\sigma _2}} )}
\end{equation}
holds for some $r > 0$.
\end{enumerate}
\end{prop}

\par

\subsection{Broader family of modulation spaces}\label{subsec1.3}

\par

\par

%As announced in the introduction we consider in Section  \ref{sec2}
%mapping properties for pseudo-differential operators when acting on
%a broader class of modulation spaces, which are defined by imposing certain
%types of translation invariant solid BF-space norms on the short-time
%Fourier transforms. (Cf. \cite{Fe4,Fe6,Fe8,FG1,FG2}.)
%
%\par

\begin{defn}\label{bfspaces1}
Let $\mascB $ be a Banach space of measurable functions on $\rr d$,
and let $v \in\mascP _E(\rr d)$.
Then $\mascB$ is called a \emph{translation invariant
Banach Function space on $\rr d$} (with respect to $v$), or \emph{invariant
BF space on $\rr d$}, if there is a constant $C$ such
that the following conditions are fulfilled:
\begin{enumerate}
\item if $x\in \rr d$ and $f\in \mascB$, then $f(\cdo -x)\in
\mascB$, and 
\begin{equation}\label{translmultprop1}
\nm {f(\cdo -x)}{\mascB}\le Cv(x)\nm {f}{\mascB}\text ;
\end{equation}

\vrum

\item if  $f,g\in L^1_{loc}(\rr d)$ satisfy $g\in \mascB$ and $|f|
\le |g|$, then $f\in \mascB$ and
$$
\nm f{\mascB}\le C\nm g{\mascB}\text ;
$$
\vrum

\item Minkowski's inequality holds true, i.{\,}e.
\begin{equation}\label{Eq:MinkIneq}
\nm {f*\fy}{\mascB}\lesssim \nm {f}{\mascB}\nm \fy{L^1_{(v)}},
\qquad f\in \mascB ,\ \fy \in L^1_{(v)}(\rr d).
\end{equation}
\end{enumerate}
\end{defn}

\par

If $v$ belongs to $\mascP _{E,s}(\rr d)$
($\mascP _{E,s}^0(\rr d)$), then $\mascB$ in Definition \ref{bfspaces1}
is called an invariant BF-space of Roumieu type (Beurling type) of order $s$.

\par

It follows from (2) in Definition \ref{bfspaces1} that if $f\in
\mascB$ and $h\in L^\infty$, then $f\cdot h\in \mascB$, and
\begin{equation}\label{multprop}
\nm {f\cdot h}{\mascB}\le C\nm f{\mascB}\nm h{L^\infty}.
\end{equation}
In Definition \ref{bfspaces1}, condition (2)  means that a
translation invariant BF-space is a solid BF-space in the sense of
(A.3) in \cite{Fe6}. 
%The space $\mascB$ in Definition \ref{bfspaces1} is called an
%\emph{invariant BF-space} (with respect to $v$) if $r=1$, and
%Minkowski's inequality holds true, i.{\,}e.
%%%
%\begin{equation}\label{Eq:MinkIneq}
%\nm {f*\fy}{\mascB}\lesssim \nm {f}{\mascB}\nm \fy{L^1_{(v)}},
%\qquad f\in \mascB ,\ \fy \in C_0^\infty (\rr d).
%\end{equation}
%%
%for some constant $C$ which is independent of
%$f\in \mascB$ and $\fy \in C_0^\infty (\rr d)$.

\par

\begin{example}\label{Lpqbfspaces}
Assume that $p,q\in [1,\infty ]$, and let $L^{p,q}_1(\rr {2d})$ be the
set of all $f\in L^1_{loc}(\rr {2d})$ such that
$$
\nm  f{L^{p,q}_1} \equiv \Big ( \int \Big ( \int |f(x,\xi )|^p\, dx\Big
)^{q/p}\, d\xi \Big )^{1/q}
$$
if finite.
%Also let $L^{p,q}_2(\rr {2d})$ be the set of all $f\in
%L^1_{loc}(\rr {2d})$ such that
%$$
%\nm  f{L^{p,q}_2} \equiv \Big ( \int \Big ( \int |f(x,\xi )|^q\, d\xi
%\Big )^{p/q}\, dx \Big )^{1/p}
%$$
%is finite.
Then it follows that $L^{p,q}_1$ %and $L^{p,q}_2$ are
is translation invariant BF-spaces with respect to $v=1$.
\end{example}

\par
%
%For translation invariant BF-spaces we make the
%following observation.
%
%\par
%
%\begin{prop}\label{p1.4BFA}
%Assume that $v\in\mascP _E(\rr {d})$, and that $\mascB$ is an
%invariant BF-space with respect to $v$ such that \eqref{Eq:MinkIneq}
%holds true. Then the
%convolution mapping $(\fy ,f)\mapsto \fy *f$ from $C_0^\infty (\rr
%d)\times \mascB$ to $\mascB$ extends uniquely to a continuous
%mapping from
%$L^1_{(v )}(\rr d)\times \mascB$ to $\mascB$, and \eqref{Eq:MinkIneq}
%holds true for any $f\in \mascB$ and $\fy \in L^1_{(v)}(\rr d)$.
%\end{prop}
%
%\par
%
%The result is a straight-forward consequence of the fact that $C_0^\infty$
%is dense in $L^1_{(v)}$.
%
%\par

We refer to \cite {Fe4,FG1,FG2,FG4,GaSa,Gc2,RSTT,To20}
for more facts about modulation spaces.
Next we consider the extended class of modulation spaces which we are interested
in.

\par

\begin{defn}\label{bfspaces2}
Assume that $\mascB$ is a translation
invariant QBF-space on $\rr {2d}$, $\omega \in\mascP _E(\rr {2d})$,
and that $\phi \in
\Sigma _1(\rr d)\setminus 0$. Then the set $M(\omega ,\mascB )$
consists of all $f\in \Sigma _1'(\rr d)$ such that
$$
\nm f{M(\omega ,\mascB )}
\equiv \nm {V_\phi f\, \omega }{\mascB}
$$
is finite.
%%% If $\omega =1$, then the notation $M(\mascB )$ is used
%%% instead of $M(\omega ,\mascB )$.
\end{defn}

\par

Obviously, we have
$
M^{p,q}_{(\omega )}(\rr d)=M(\omega ,\mascB )$
when  $\mascB =L^{p,q}_1(\rr {2d})$ (cf. Example \ref{Lpqbfspaces}).
It follows that many properties which are valid for the classical modulation
spaces also hold for the spaces of the form $M(\omega ,\mascB )$.
%For example we have the following proposition, which shows that
%the definition of $M(\omega ,\mascB )$ is independent of the
%choice of $\phi$ when $\mascB$ is a Banach space. We omit the proof
%since it follows by similar arguments as
%in the proof of Proposition 11.3.2 in \cite{Gc2}.

\par

We notice that $M(\omega ,\mascB )$ is independent of the choice
of $\phi$ in Definition \ref{bfspaces2} cf.\, \cite{To25}.
Furthermore, $M(\omega ,\mascB )$ is a Banach space in view of \cite{PfTo}.

\par

\subsection{Pseudo-differential operators}

\par

Next we recall some facts on pseudo-differential operators. Let
$A\in \GL (d,\mathbf R)$ be fixed and
let $a\in \Sigma _1(\rr {2d})$. Then the pseudo-differential
operator $\op _A(a)$ is the linear and continuous operator on $\Sigma _1(\rr d)$,
defined by the formula
\begin{multline}\label{e0.5}
(\op _A(a)f)(x)
\\[1ex]
=
(2\pi  ) ^{-d}\iint a(x-A(x-y),\xi )f(y)e^{i\scal {x-y}\xi }\,
dyd\xi .
\end{multline}
The definition of $\op _A(a)$ extends to
any $a\in \Sigma _1'(\rr {2d})$, and then $\op _A(a)$ is continuous from
$\Sigma _1(\rr d)$ to $\Sigma _1'(\rr d)$. Moreover, for every fixed
$A\in \GL (d,\mathbf R)$, it follows that there is a one to
one correspondence between such operators and pseudo-differential
operators of the form $\op _A(a)$. (See e.{\,}g. \cite {Ho1}.)
If $A=2^{-1}I$, where $I\in \GL (d,\mathbf R)$ is the identity matrix, then
$\op _A(a)$ is equal to the Weyl operator $\op ^w(a)$
of $a$. If instead $A=0$, then the standard (Kohn-Nirenberg)
representation $\op (a)$ is obtained.

\par

If $a_1,a_2\in \Sigma _1'(\rr {2d})$ and $A_1,A_2\in
 \GL (d,\mathbf R)$, then
\begin{equation}\label{pseudorelation}
\op _{A_1}(a_1)=\op _{A_2}(a_2) \quad \Leftrightarrow \quad a_2(x,\xi
)=e^{i\scal {(A_1-A_2)D_\xi}{D_x}}a(x,\xi ).
\end{equation}
(Cf. \cite{Ho1}.)

\par

%The following special case of \cite[Theorem 3.1]{To24} is important when discussing
%continuity of pseudo-differential operators when acting on quasi-Banach modulation
%spaces.
%
%\par
%
%\begin{prop}\label{Prop:OpCont}
%Let $\omega _1,\omega _2
%\in \mathscr P_{E}(\rr {2d})$ and $\omega _0\in \mathscr P_{E}(\rr {2d}\oplus \rr {2d})$
%be such that
%%%
%\begin{equation}\label{Eq:WeightfracCond1}
%\frac {\omega _2(x,\xi  )}{\omega _1
%(y,\eta )} \lesssim \omega _0( x,\eta ,\xi -\eta ,y-x ).
%\end{equation}
%%%
%Also let $\mabfp \in (0,\infty]^{2d}$,
%$E$ be a phase split basis of $\rr {2d}$  and let
%$a\in M^{\infty ,1}_{(\omega _0)}(\rr {2d})$. Then
%$\op _0(a)$ is continuous from $M(L^{\mabfp ,E},
%\omega _1)$ to $M(L^{\mabfp ,E},\omega _2)$.
%\end{prop}
%
%\par
%
%In the next section we discuss continuity for pseudo-differential
%operators with symbols in the following definition. (See also the introduction.)

\par

%%%%%
\subsection{Symbol classes}
%%%%%

\par

Next we introduce function spaces related to symbol classes
of the pseudo-differential operators. These functions should obey various
conditions of the form
\begin{align}
|\partial _x^\alpha \partial _\xi ^\beta a(x,\xi )|
&\lesssim
h ^{|\alpha +\beta |}\alpha !^\sigma \beta !^s  \omega (x,\xi ),
\label{Eq:symbols2}
\end{align}
for functions on $\rr {d_1+d_2}$. For this reason we consider
semi-norms of the form
\begin{equation}\label{Eq:GammaomegaNorm}
\nm a{\Gamma _{(\omega )}^{\sigma ,s;h}}
\equiv \sup _{(\alpha ,\beta) \in \nn {d_1+d_2}}
\left (
\sup _{(x,\xi) \in \rr {d_1+d_2}} \left (
\frac {|\partial _x^\alpha \partial _\xi ^\beta a(x,\xi )|}{
	h ^{|\alpha +\beta |}\alpha !^\sigma \beta !^s \omega (x,\xi )}
\right ) \right ) ,
\end{equation}
indexed by $h>0$, 

\par

\begin{defn}\label{Def:GammaSymb2}
Let $s$, $\sigma$ and $h$ be positive constants,
let $\omega$ be a weight on $\rr {d_1+d_2}$, and let
$$
\omega _r(x,\xi )\equiv e^{r(|x|^{\frac 1s} + |\xi |^{\frac 1\sigma })}.
$$ 
\begin{enumerate}
\item The set $\Gamma _{(\omega )}^{\sigma ,s;h}  (\rr {d_1+d_2})$
consists of
all $a \in C^\infty(\rr {d_1+d_2})$ such that
$\nm a{\Gamma _{(\omega )}^{\sigma ,s;h}}$ in
\eqref{Eq:GammaomegaNorm} is finite.
The set $\Gamma _0^{\sigma ,s;h}  (\rr {d_1+d_2})$ consists of
all $a \in C^\infty(\rr {d_1+d_2})$ such that
$\nm a{\Gamma _{(\omega_r )}^{\sigma ,s;h}}$ is finite
for every $r>0$, and the topology is the projective limit topology of
$\Gamma _{(\omega _r)}^{\sigma ,s;h}  (\rr {d_1+d_2})$ with respect to $r>0$;

\vrum

\item The sets $\Gamma _{(\omega )}^{\sigma ,s}  (\rr {d_1+d_2})$ and
$\Gamma _{(\omega )}^{\sigma ,s;0}  (\rr {d_1+d_2})$ are given by
\begin{align*}
\Gamma _{(\omega )}^{\sigma ,s}  (\rr {d_1+d_2})
&\equiv
\bigcup _{h>0}\Gamma _{(\omega )}^{\sigma ,s;h}  (\rr {d_1+d_2})
\intertext{and}
\Gamma _{(\omega )}^{\sigma ,s;0}  (\rr {d_1+d_2})
&\equiv
\bigcap _{h>0}\Gamma _{(\omega )}^{\sigma ,s;h}  (\rr {d_1+d_2}),
\end{align*}
and their topologies are the inductive respective the projective topologies
of $\Gamma _{(\omega )}^{\sigma ,s;h}  (\rr {d_1+d_2})$ with respect to $h>0$.	
\end{enumerate}	
\end{defn}

\par

%\begin{rem}\label{Remark:Gammaclasses}
%We have
%$$
%\mascP \subseteq \mascP _{E,s_1}^0\subseteq \mascP _{E,s_1}
%\subseteq \mascP _{E,s_2}^0, \qquad s_2<s_1.
%$$
%Hence, despite that $\Gamma ^{(\omega _0)}_{0,s}\subseteq
%\Gamma ^{(\omega _0)}_s\subseteq S^{(\omega _0)}$ holds for every
%$\omega _0$,  we have
%%%
%\begin{align*}
%\Gamma ^{(\omega )}_{0,s} &\nsubseteq
%\bigcup _{\omega _0\in \mascP} S^{(\omega _0)}
%\intertext{for some $\omega \in \mascP _{E,s}^0$, and}
%\Gamma ^{(\omega )}_{0,s} &\nsubseteq
%\bigcup _{\omega _0\in \mascP _{E,s}^0} \Gamma ^{(\omega _0)}_{s}
%\end{align*}
%%%
%for some $\omega \in \mascP _{E,s}$.
%\end{rem}
%
%\par
\par

The following result is a straight-forward consequence of 
\cite[Proposition 2.4]{AbCaTo} and the definitions.

\par

\begin{prop}\label{SymbClassModSpace}
Let $R>0$, $q\in (0,\infty ]$, $s,\sigma >0$ be such that $s+\sigma \ge 1$
and $(s,\sigma )\neq (\frac 12,\frac 12)$, $\phi \in \Sigma _{s,\sigma}^{\sigma ,s}
(\rr {2d})\setminus 0$,
$\omega \in \mascP _{s,\sigma}(\mathbf{R}^{2d}),$ and let
$$
\omega _R(x,\xi, \eta, y) = \omega(x,\xi) e^{-R(|y|^{\frac 1s} + |\eta|^{\frac 1\sigma})}.
$$
Then
\begin{align}\label{iden}
\begin{split}
\Gamma^{\sigma ,s}_{(\omega)}(\rr {2d})
&=
\bigcup _{R>0}\sets {a\in (\Sigma _{s,\sigma}^{\sigma ,s})'
	(\rr {2d})}{\nm {\omega
		_R^{-1}V_\phi a}{L^{\infty ,q}} <\infty },
\\[1ex]
\Gamma^{\sigma ,s;0}_{(\omega)}(\rr {2d})
&=
\bigcap _{R>0}\sets {a\in (\Sigma _{s,\sigma}^{\sigma ,s})'
	(\rr {2d})}{\nm {\omega
		_R^{-1}V_\phi a}{L^{\infty ,q}} <\infty }.
\end{split}
\end{align}
\end{prop}

\par

The following lemma is a consequence of \cite[Theorem 3.6]{AbCaTo}.
 
\par

\begin{lemma}\label{Somega}
Let $s,\sigma>0$ such that $s+\sigma\ge 1$ $\omega \in \mascP _{s,\sigma}(\rr {2d})$,
$A_1,A_2\in \GL (d,\mathbf R )$, and
that $a_1,a_2\in (\Sigma _{s,\sigma}^{\sigma,s})'(\rr {2d})$ are such that
$\op _{A_1}(a_1)=\op _{A_2}(a_2)$. Then
\begin{alignat*}{3}
a_1 &\in \Gamma _{(\omega )}^{\sigma,s;0}(\rr {2d})&\qquad &\Leftrightarrow &\qquad a_2
&\in \Gamma_{(\omega )}^{\sigma,s;0}(\rr {2d}) 
%\intertext{and}
%a_1 &\in \Gamma _{(\omega )}^{\sigma,s}(\rr {2d})&\qquad &\Leftrightarrow &\qquad a_2
%&\in \Gamma _{(\omega )}^{\sigma,s}(\rr {2d}) .
\end{alignat*}
and similarly for $\Gamma _{(\omega )}^{\sigma,s}(\rr {2d})$ in
place of $\Gamma _{(\omega )}^{\sigma,s;0}(\rr {2d})$.
\end{lemma}

\par
%%%%%%%%%%%%%%%%%%%%%%%%%%%%%%%%%%
\section{Continuity for pseudo-differential operators
with symbols of infinite order}\label{sec2}
%%%%%%%%%%%%%%%%%%%%%%%%%%%%%%%%%%

\par

In this section we discuss continuity for operators in
$\op (\Gamma ^{\sigma,s}_{(\omega_0)})$ and
$\op (\Gamma ^{\sigma,s;0}_{(\omega _0)})$
when acting on a general class of modulation spaces.
In Theorem \ref{p3.2} continuity is treated where the symbols belong to
$\Gamma ^{\sigma,s}_{(\omega _0)}$ and in Theorem \ref{p3.2B}
 continuity is treated where the symbols belong to
$\Gamma ^{\sigma,s;0}_{(\omega _0)}$.
%
%below it is proved that if $\omega ,\omega
%_0\in\mascP _{s,\sigma}^0$, $A\in \GL (d,\mathbf R)$ and
%$a\in \Gamma ^{\sigma,s}_{(\omega _0)}$, then
%$\op _A(a)$ is continuous from $M(\omega _0\omega ,\mascB )$ to
%$M(\omega ,\mascB )$.
%
This gives an analogy to \cite[Theorem 3.2]{To14}
in the framework of operator theory and Gelfand-Shilov classes.

\par

Our main result is stated as follows.

\par

\begin{thm}\label{p3.2}
Let $A\in \GL (d,\mathbf R)$, $s,\sigma\ge 1$,
$\omega ,\omega _0\in\mascP _{s,\sigma}^0(\rr {2d})$,
$a\in \Gamma _{(\omega _0)}^{\sigma,s}(\rr {2d})$, and that $\mascB$
is an invariant BF-space on $\rr {2d}$. Then
$\op _A(a)$ is continuous from $M(\omega _0\omega ,\mascB )$
to $M(\omega ,\mascB )$.
\end{thm}

\par

We need some preparations for the proof, and start with the following remark.

%Recalling Minkowski's
%inequality in a somewhat general form. Assume that
%$d\mu$ is a positive measure, and that $f\in L^1(d\mu ;\mascB )$
%for some Banach space $\mascB$. Then Minkowski's inequality asserts that
%$$
%\Big \Vert \int f(x)\, d\mu (x)\Big \Vert _{\mascB} \le \int  \nm {f(x)}{\mascB}\, d\mu (x).
%$$
%
%\par
%
%We also need the following remark and some lemmas.

\par

\begin{rem}\label{rem:adjUniq}
Let $s,\sigma>0$ such that $s+\sigma \geq 1$.
If $a\in (\Sigma _{s,\sigma}^{\sigma,s})'(\rr {2d})$,
then there is a unique $b\in (\Sigma _{s,\sigma}^{\sigma,s})'(\rr {2d})$ such that
$\op (a)^*=\op (b)$, where $b(x,\xi )= e^{i\scal {D_\xi}{D_x}}\overline {a(x,\xi )}$
in view of \cite[Theorem 18.1.7]{Ho1}. Furthermore, by the latter equality and
\cite[Theorem 4.1]{CaTo} it follows that
$$
a\in \Gamma _{(\omega )}^{\sigma,s}(\rr {2d})
\quad \Leftrightarrow \quad
b\in \Gamma _{(\omega )}^{\sigma,s}(\rr {2d}).
$$
\end{rem}

\par

\par

\begin{lemma}\label{lem:equivfun}
Suppose $s,\sigma\geq 1$,
$\omega \in \mascP _E(\rr {d_0})$ and that $f\in C^\infty
(\rr {d+d_0})$ satisfies
\begin{equation}\label{eq:anGeShiEst}
|\partial ^\alpha f(x,y)|\lesssim h^{|\alpha |}\alpha !^\sigma 
e^{-r|x|^{\frac 1s}}\omega (y),	\alpha \in \nn {d+d_0}
\end{equation}
for some $h>0$ and $r>0$. Then there are $f_0\in C^\infty (\rr {d+d_0})$
and $\psi \in \maclS _s^\sigma(\rr d)$ such that \eqref{eq:anGeShiEst} holds
with $f_0$ in place of $f$ for some for some $h>0$ and $r>0$, and
$f(x,y)= f_0(x,y)\psi (x)$.		
\end{lemma}

\par

\begin{proof}
By Proposition \ref{Prop:EquivWeights}, there is a submultiplicative weight
$v_0\in \mascP _{E,s}(\rr d)\cap C^\infty (\rr d)$ such that
\begin{align}
v_0(x)&\asymp e^{\frac r2|x |^{\frac 1s}}\label{eq:v0Est}
\intertext{and}
|\partial ^\alpha v_0(x)| &\lesssim h^{|\alpha |}\alpha !^\sigma 
v_0(x),\qquad \alpha \in \nn d
\label{eq:vEst}
\intertext{for some $h,r>0$.	
Since $s,\sigma\ge 1$,
a straight-forward application of Fa{\`a} di Bruno's formula,
for the composed function $\psi(x)=g(v_0(x))$, where $g(t)=\frac 1t$,
on \eqref{eq:vEst} gives}
\left | \partial ^\alpha \left (\frac 1{v_0(x)}\right ) \right | &\lesssim
h^{|\alpha |}\alpha !^\sigma\cdot \frac 1{v_0(x)},\qquad \alpha \in \nn d
\tag*{(\ref{eq:vEst})$'$}
\end{align}
for some $h>0$. It follows from \eqref{eq:v0Est} and \eqref{eq:vEst}$'$
that if $\psi =1/{v_0}$, then $\psi \in \maclS _s^\sigma(\rr d)$.
Furthermore, if $f_0(x,y)=f(x,y)v_0(x)$,
then an application of Leibnitz formula we get
\begin{multline*}
|\partial ^{\alpha}_x\partial ^{\alpha _0}_yf_0(x,y)|\lesssim
\sum _{\gamma \le \alpha }\binom{\alpha}{\gamma} |\partial ^\delta _x
\partial ^{\alpha _0}_yf(x,y)|
\, |\partial ^{\alpha -\delta }v_0(x)|
\\[1ex]
\lesssim
h^{|\alpha |+|\alpha _0|} \sum _{\gamma \le \alpha }
\binom{\alpha}{\gamma}
(\gamma !\alpha _0!)^\sigma e^{-r|x|^{\frac 1s}}
\omega (y)(\alpha -\gamma )!^\sigma v_0(x)
\\[1ex]
\lesssim
(2h)^{|\alpha |+|\alpha _0|}(\alpha !\alpha _0!)^\sigma
e^{-r|x|^{\frac 1s}}v_0(x)\omega (y)
\\[1ex]
\lesssim
(2h)^{|\alpha |+|\alpha _0|}(\alpha !\alpha _0!)^\sigma
e^{-\frac r2|x|^{\frac 1s}}\omega (y)
\end{multline*}
for some $h>0$, which gives the desired estimate on $f_0$,
since it is clear that $f(x,y)= f_0(x,y)\psi (x)$.
\end{proof}

\par

\begin{lemma}\label{Lemma:PrepReThm3.2}
Let $s,\sigma \ge 1$, $\omega \in
\mascP _{s,\sigma}^0(\rr {2d})$, $v_1\in
\mascP _{s}^0(\rr {d})$ and $v_2\in \mascP _{\sigma}^0(\rr d)$ be such that
$v_1$ and $v_2$ are submultiplicative, $\omega \in \Gamma _{(\omega )}
^{\sigma,s}(\rr {2d})$ is $v_1\otimes v_2$-moderate.
Also let $a\in \Gamma _{(\omega )}^{\sigma,s}(\rr {2d})$, 
$f\in \maclS _s^\sigma(\rr d)$, $\phi \in \Sigma _s^\sigma(\rr d)$,
$\phi_2=\phi v_1$,
If 
\begin{align}
\Phi (x,\xi ,z,\zeta ) &= \frac {a(x+z,\xi +\zeta )}
{\omega (x,\xi)v_1(z)v_2(\zeta )}\label{eq:Phidef}
\intertext{and}
H(x,\xi ,y) &= \iint \Phi (x,\xi ,z,\zeta )\phi _2(z)v_2(\zeta )
e^{i\scal {y-x-z}{\zeta}}\, dzd\zeta .\label{eq:Hidentity}
\end{align}
Then
\begin{equation}\label{eq:stftpseudoform}
V_\phi (\op (a)f)(x,\xi ) = (2\pi )^{-d} (f,e^{i\scal \cdo \xi
}H(x,\xi ,\cdo ))\omega (x,\xi
).
\end{equation}
Furthermore the following is true:
\begin{enumerate}
\item $H\in C^\infty (\rr {3d})$ and satisfies
\begin{equation}\label{eq:DerHEst}
|\partial _y^\alpha H(x,\xi ,y)| \lesssim h_0^{|\alpha |}
\alpha!^\sigma e^{-r_0|x-y|^{\frac 1s}},
\end{equation}
for every $\alpha \in \rr d$ and some $h_0,r_0>0$;
\vrum

\item there are functions $H_0\in C^\infty (\rr {3d})$
and $\phi _0\in \maclS _s^\sigma(\rr d)$
such that
\begin{equation}\label{eq:HProd}
H(x,\xi ,y) = H_0(x,\xi ,y)\phi _0(y-x),
\end{equation}
and such that \eqref{eq:DerHEst} holds for some $h_0,r_0>0$,
with $H_0$ in place of $H$.		
\end{enumerate}	
\end{lemma}

\par

Lemma \ref{Lemma:PrepReThm3.2} follows by similar arguments as in
\cite{To25}. In order to be self contained we give a different proof.

\par

\begin{proof}
%When proving the first part, we shall use some ideas in the proof of
%Lemma \ref{lem:STFTequivFormula}.
By straight-forward computations we get
\begin{equation}\label{eq:storformel}
V_\phi (\op (a)f)(x,\xi ) 
%	=(a(\cdot ,D)f, 	\phi (\cdo -x)\, e^{i\scal \cdot \xi})
%	\\[1ex]
%	= (f, \op (a)^*(\phi (\cdo -x)\, e^{i\scal \cdot \xi}))
%	\\[1ex]
= (2\pi )^{-d} (f,e^{i\scal \cdo \xi }H_1(x,\xi ,\cdo
))\omega (x,\xi ),
\end{equation}
where
\begin{multline*}
H_1(x,\xi ,y) = (2\pi )^{d}e^{-i\scal y\xi}(\op (a)^*(\phi (\cdo -x)
\, e^{i\scal \cdot \xi}))(y)/\omega (x,\xi )
\\[1ex]
= \iint \frac {a(z,\zeta )}{\omega (x,\xi )}\phi (z-x)
e^{i\scal {y-z}{\zeta -\xi}}\, dzd\zeta 
\\[1ex]
= \iint \Phi (x,\xi ,z-x,\zeta -\xi )\phi _2(z-x)v_2(\zeta
-\xi )e^{i\scal {y-z}{\zeta -\xi}}\, dzd\zeta .
\end{multline*}
If $z-x$ and $\zeta -\xi$ are taken as new variables of integrations,
it follows that the right-hand side is the same as \eqref{eq:Hidentity}.
Hence \eqref{eq:stftpseudoform} holds. This
gives the first part of the lemma.

\par

The smoothness of  $H$ is a consequence of the uniqueness of the adjoint
(cf. Remark \ref{rem:adjUniq}) and \cite[Lemma 2.7]{To25}.

\par

To show that \eqref{eq:DerHEst} holds, let
$$
\Phi _0(x,\xi ,z,\zeta ) = \Phi (x,\xi ,z,\zeta )\phi _2(z),
$$
where $\Phi$ defined as in \eqref{eq:Phidef},
and let $\Psi =\mascF _3\Phi_0$, where $\mascF _3\Phi_0$ is the partial
Fourier transform of $\Phi _0(x,\xi ,z,\zeta )$ with respect to the $z$ variable.
Then it follows from the assumptions and \eqref{eq:vEst}$'$ that
\begin{multline*}
|\partial _z^\alpha \Phi _0(x,\xi ,z,\zeta )|
%=\left|
%\sum_{\gamma\leq \alpha}\binom{\alpha}{\gamma}
%\partial_z^{\gamma}\prn{\frac{a(x+z,\xi+\zeta)}{\omega(x,\xi)v_1(z)v_2(\zeta)}}
%\partial^{\alpha-\gamma}\phi _2(z)
%\right|
%\\
\lesssim \sum_{\gamma\leq \alpha}\binom{\alpha}{\gamma}
\sum _{\lambda\leq \gamma}\binom{\gamma}{\lambda}
\frac{\left|\partial_z^{\gamma-\lambda}a(x+z,\xi+\zeta)\right|}{\omega(x,\xi)v_2(\zeta)}
\\
\times \partial^{\lambda}\prn{\frac{1}{v_1(z)}}
h^{|\alpha-\gamma|}(\alpha-\gamma)!^\sigma e^{-r|z|^{\frac 1s}}
%\\[1ex]
%\lesssim \sum_{\gamma\leq \alpha}\binom{\alpha}{\gamma}
%\sum _{\lambda\leq \gamma}\binom{\gamma}{\lambda}
%h^{|\alpha-\lambda|}(\alpha-\gamma)!^\sigma(\gamma-\lambda)!^\sigma
%v_1(z)\cdo \partial^{\lambda}\prn{\frac{1}{v_1(z)}}
%e^{-r|z|^{\frac 1s}}
\\[1ex]
\lesssim \sum_{\gamma\leq \alpha}\binom{\alpha}{\gamma}
\sum _{\lambda\leq \gamma}\binom{\gamma}{\lambda}
h^{|\alpha|}(\alpha-\gamma)!^\sigma(\gamma-\lambda)!^\sigma \lambda!^\sigma
e^{-r_0|z|^{\frac 1s}}
\\[1ex]
\lesssim  h^{|\alpha|} \alpha!^\sigma
\sum_{\gamma\leq \alpha}\binom{\alpha}{\gamma}
\sum _{\lambda\leq \gamma}\binom{\gamma}{\lambda}
\prn{\frac{(\alpha-\gamma)!\gamma!}{\alpha!}}^\sigma
\prn{\frac{(\gamma-\lambda)!\lambda!}{\gamma!}}^\sigma
e^{-r_0|z|^{\frac 1s}}
%\\[1ex]
%\lesssim  h^{|\alpha|} \alpha!^\sigma e^{-r_0|z|^{\frac 1s}}
%\sum_{\gamma\leq \alpha}2^{|\alpha|}
%\sum _{\lambda\leq \gamma}2^{|\gamma|}
\\[1ex]
\lesssim  (4h)^{|\alpha|} \alpha!^\sigma e^{-r|z|^{\frac 1s}}
\sum_{\gamma\leq \alpha}1\cdo \sum _{\lambda\leq \gamma}1.
\end{multline*}
Since $\sum _{\lambda\leq \gamma}1 \lesssim 2^{|\gamma|}$,
we get
\begin{equation}\label{eq:Phi_0Est}
|\partial _z^\alpha \Phi _0(x,\xi ,z,\zeta )|
\leq C (16h)^{|\alpha|}\alpha!^\sigma e^{-r_0|z|^{\frac 1s}}
\leq C h_0^{|\alpha|}\alpha!^\sigma e^{-r_0|z|^{\frac 1s}}
\end{equation}
for some $C,h_0,r_0>0$. Then $z\mapsto \Phi _0(x,\xi ,z,\zeta )$
is an element in $\maclS _s^\sigma(\rr d)$.
Moreover, $\set{\Phi _0(x,\xi ,z,\zeta )}_{z\in\rr d}$ is a bounded set in
$\Gamma _{(1)}^{\sigma,s}(\rr d\times\rr {2d})$.
Indeed, for a fixed $z_0\in \rr{d}$, then an application of Leibnitz formula,
Fa{\`a} di Bruno's formula, Proposition \ref{Prop:EquivWeights} and \eqref{eq:vEst}$'$, give
\begin{multline*}
\abs{\partial_x^\alpha\partial_\xi^\beta\partial_\zeta^\gamma\Phi _0(x,\xi ,z_0,\zeta )}
\leq \sum
\binom \alpha{\alpha_1}\binom \beta{\beta_1}\binom \gamma{\gamma_1}
\partial_x^{\alpha_1}\partial_\xi^{\beta_1}\prn{\frac 1{\omega(x,\xi)}}
\\%[1ex]
\times\partial_\zeta^{\gamma_1}\prn{\frac 1{v_2(\zeta)}}\abs{\partial_x^{\alpha-\alpha_1}
	\partial_\xi^{\beta-\beta_1}\partial_\zeta^{\gamma-\gamma_1}
	a(x+z_0,\xi+\zeta)}\cdo \frac{|\phi(z_0)|}{v_1(z_0)}
\\[1ex]
\lesssim
\sum
\binom \alpha{\alpha_1}\binom \beta{\beta_1}\binom \gamma{\gamma_1}
h^{|\alpha_1+\beta_1+\gamma_1|}\alpha_1!^\sigma(\beta_1!\gamma_1!)^s
\\%[1ex]
\times\prn{\frac 1{\omega(x,\xi)v_1(z_0)v_2(\zeta)}}
\abs{\partial_x^{\alpha-\alpha_1}
	\partial_\xi^{\beta-\beta_1}\partial_\zeta^{\gamma-\gamma_1}
	a(x+z_0,\xi+\zeta)}
\\[1ex]
\lesssim
h^{|\alpha+\beta+\gamma|}\sum
\binom \alpha{\alpha_1}\binom \beta{\beta_1}\binom \gamma{\gamma_1}
\prn{(\alpha-\alpha_1)!\alpha_1!}^\sigma
\prn{(\beta-\beta_1)!\beta_1!}^s\prn{(\gamma-\gamma_1)!\gamma_1!}^s
%\\[1ex]
%\lesssim
%(2h)^{|\alpha+\beta+\gamma|}\alpha!^\sigma(\beta!\gamma!)^s
%\sum_{\substack{\alpha_1\leq \alpha\\
%		\beta_1\leq\beta\\\gamma_1\leq\gamma}}
%\prn{\frac{(\alpha-\alpha_1)!\alpha_1!}{\alpha!}}^\sigma
%\prn{\frac{(\beta-\beta_1)!\beta_1!}{\beta!}}^s
%\prn{\frac{(\gamma-\gamma_1)!\gamma_1!}{\gamma!}}^s
\\[1ex]
\lesssim
(4h)^{|\alpha+\beta+\gamma|}\alpha!^\sigma(\beta!\gamma!)^s,   	
\end{multline*}
where all the summations above are taken over all 
$\alpha_1\leq \alpha, \beta_1\leq\beta$ and $\gamma_1\leq\gamma$.
In view of Proposition \ref{propBroadGSSpaceChar}
and \eqref{eq:Phi_0Est} we have 
$$
|\partial _\eta^\alpha \Psi (x,\xi ,\eta ,\zeta )|
\lesssim h_0^{|\alpha |}\alpha !^se^{-r_0|\eta |^{\frac 1\sigma}},
$$
for some $h_0,r_0>0$. Hence
$$
|\partial _\eta^\alpha (\Psi (x,\xi ,\zeta ,\zeta )v_2(\zeta ))|
\lesssim h_0^{|\alpha |}\alpha !^se^{-r_0|\zeta |^{\frac 1\sigma}}
$$
for some $h_0,r_0>0$.

\par

By letting $H_2(x,\xi ,\cdo )$ be the inverse partial Fourier transform of
$\Psi (x,\xi ,\zeta ,\zeta )v_2(\zeta )$ with respect to the $\zeta$ variable,
it follows that
\begin{equation}\label{eq:H2Est}
|\partial _y^\alpha H_2(x,\xi ,y)|
\lesssim h_0^{|\alpha |}\alpha !^\sigma e^{-r_0|y|^{\frac 1s}}
\end{equation}
for some $h_0,r_0>0$. The assertion (1) now follows from the latter estimate
and the fact that $H(x,\xi ,y)= H_2(x,\xi ,x-y)$.

\par

In order to prove (2) we notice that \eqref{eq:H2Est} shows that
$y\mapsto H_2(x,\xi ,y)$ is an element in $\maclS _s^\sigma(\rr d)$ with values
in $\Gamma ^{(1)}_{s,\sigma}(\rr {2d})$. It follows by Lemma \ref{lem:equivfun}
that there exist $H_3\in C^\infty (\rr {3d})$ and $\phi _0\in \maclS _s^\sigma(\rr d)$
such that \eqref{eq:H2Est} holds for some $h_0,r_0>0$ with $H_3$ in place of $H_2$,
and
$$
H_2(x,\xi ,y)= H_3(x,\xi ,y)\phi _0(-y).
$$
This is the same as (2), and the result follows.	
\end{proof}

\par

\begin{proof}[Proof of Theorem \ref{p3.2}]
There is no restriction if we assume that $A=0$.
Let $G=\op (a)f$. In view of Lemma \ref{Lemma:PrepReThm3.2} we have
\begin{multline*}
V_\phi G(x,\xi )
=
(2\pi )^{-\frac d2}\mascF ((f\cdot \overline {\phi _0(\cdo -x)}) \cdot H_0(x,\xi ,\cdo ))(\xi )
\omega (x,\xi )
%\\[1ex]
%=
%(2\pi )^{-d}\mascF ((f\cdot \overline {\phi _0(\cdo -x)})) *(\mascF  (H_0(x,\xi ,\cdo )))(\xi )
%\omega (x,\xi )
\\[1ex]
=
(2\pi )^{-d} (V_{\phi _0}f)(x,\cdo ) * (\mascF  (H_0(x,\xi ,\cdo )))(\xi )
\omega (x,\xi ).
\end{multline*}
Since $\omega$ and $\omega _0$ belong to $\mascP _{s,\sigma}^0(\rr {2d})$,
then for every $r_0>0$ and $x,\xi,\eta\in \rr d$ we have
$$
\omega (x,\xi )\omega _0(x,\xi )
\lesssim
\omega (x,\eta )\omega _0(x,\eta ) e^{\frac {r_0}2|\xi-\eta|^{\frac 1\sigma}},
$$

this inequality and (2) in Lemma \ref{Lemma:PrepReThm3.2} give
\begin{equation*}
|V_\phi G(x,\xi )\omega _0(x,\xi )|
\lesssim
\left(
|(V_{\phi _0}f)(x,\cdo )\omega (x,\cdo )
\omega _0(x,\cdo )| * e^{-\frac {r_0}2|\cdo |^{\frac 1\sigma}}
\right)
(\xi).
\end{equation*}

\par

In view of Definition \ref{bfspaces1},
we get for some $v\in \mascP _{\sigma}^0(\rr d)$,
\begin{multline*}
\nm G{M(\omega _0,\mascB )}\lesssim \nm {|(V_{\phi _0}f) 
\cdot \omega \cdot \omega _0|
* \delta _0 \otimes e^{-r_0|\cdo |^{\frac 1\sigma}}}{\mascB}
\\[1ex]
\le
\nm {(V_{\phi _0}f) \cdot \omega \cdot \omega _0}{\mascB}
\nm {e^{-r_0|\cdo |^{\frac 1\sigma}}v}{L^1}
\asymp \nm f{M(\omega \cdot \omega _0,\mascB )}.
\end{multline*}
This gives the result.
\end{proof}

\par

By similar arguments as in the proof of Theorem \ref{p3.2}
and Lemma \ref{Lemma:PrepReThm3.2}  we get the following.
The details are left for the reader.

\par

\begin{thm}\label{p3.2B}
Let $A\in \GL (d,\mathbf R)$, $s,\sigma\ge 1$,
$\omega ,\omega _0\in\mascP _{s,\sigma}(\rr {2d})$,
$a\in \Gamma _{(\omega _0)}^{\sigma,s;0}(\rr {2d})$, and that $\mascB$
is an invariant BF-space on $\rr {2d}$. Then
$\op _A(a)$ is continuous from $M(\omega _0\omega ,\mascB )$
to $M(\omega ,\mascB )$.
\end{thm}

\par

\begin{lemma}\label{Lemma:PrepReThm3.2B}
Let $s,\sigma \ge 1$, $\omega \in
\mascP _{s,\sigma}(\rr {2d})$, $v_1\in
\mascP _{s}(\rr {d})$ and $v_2\in \mascP _{\sigma}(\rr d)$ be such that
$v_1$ and $v_2$ are submultiplicative, $\omega \in \Gamma _{(\omega )}
^{\sigma,s;0}(\rr {2d})$ is $v_1\otimes v_2$-moderate.
Also let $a\in \Gamma _{(\omega )}^{\sigma,s;0}(\rr {2d})$, 
$f, \phi \in \Sigma _s^\sigma(\rr d)$,
$\phi_2=\phi v_1$, and let $\Phi$ and $H$ be as in Lemma
\ref{Lemma:PrepReThm3.2}.
Then \eqref{eq:stftpseudoform} and the following hold true:
\begin{enumerate}
\item $H\in C^\infty (\rr {3d})$ and satisfies \eqref{eq:DerHEst}
for every $h_0,r_0>0$;

\vrum

\item there are functions $H_0\in C^\infty (\rr {3d})$
and $\phi _0\in \Sigma _s(\rr d)$
such that \eqref{eq:HProd} holds,
and such that \eqref{eq:DerHEst} holds for every $h_0,r_0>0$,
with $H_0$ in place of $H$.
\end{enumerate}
\end{lemma}

\par

\end{document}